\documentclass[11pt]{amsart}

\usepackage{amssymb}

\newtheorem{Lemma}{Lemma}

\newtheorem{Theorem}[Lemma]{Theorem}
\newtheorem{Corollary}[Lemma]{Corollary}
\newtheorem{Proposition}[Lemma]{Proposition}
\newtheorem{Conjecture}[Lemma]{Conjecture}
\newtheorem{Question}[Lemma]{Question}

\newtheorem{Example}[Lemma]{Example}

\newtheorem*{theorem*}{Theorem}

\renewcommand{\>}{\rangle}
\newcommand{\Z}{\mathbb{Z}}

\begin{document}

\title[Two remarks on elementary theories of groups]{Two remarks on elementary theories of groups obtained by free constructions} 

\author{Eric Jaligot}
\address{Institut Fourier, 
Universit\'e de Grenoble 1, 
100 rue des maths, BP 74, 
38402 St Martin d'H\`eres cedex, 
France}

\date{\today} 

\thanks{Work partially supported be the ERC Grant 278722-GTMT}

\begin{abstract}
We give two slight generalizations of results of Poizat about elementary theories of groups obtained 
by free constructions. The first-one concerns generic types and 
the non-superstability of such groups in many cases. 
The second-one concerns the connectedness of most free products of groups without amalgamation. 
\end{abstract}

\maketitle 

First-order theories of free products of groups have been recently investigated in 
\cite{MR2776984} and \cite{Sela10}, with some transfers of arguments from free groups to 
free products of groups. We expect that some of this work transfers further to more general 
classes of groups obtained by free constructions. In this modest and short note 
we will make slight generalizations of early arguments of Poizat on first-order theories of free groups. 

By a group obtained by a free construction, we mean a free product of groups with amalgamation, 
or an $HNN$-extension. Such groups are best analized by the Bass-Serre theory of actions on 
trees, and we refer to \cite{Houcine2008} for links with Model Theory. Here our approach is 
purely algebraic and will adapt (and explain) some arguments contained in \cite{MR704088}. 

We recall that one of the main accomplishements of \cite{Sela10} is a proof of the stability 
of the free product of two stable groups. 
Conversely, we point out the following question as a possibly difficult one. 

\begin{Question}\label{Question1}
Can one have a free product of groups $G\ast H$ with a stable theory, but with the factor $G$ unstable?
\end{Question} 

\noindent
We expect that the stability result in \cite{Sela10} generalizes to certain free products 
with amalgamation and certain $HNN$-extensions, and Question \ref{Question1} also makes sense in 
these more general cases. 
When considering a stable group obtained by a free construction below, we will not assume the stability of the factors. 

In a first series of results we make a basic analysis of generic types in many groups obtained by 
free constructions, and derive the non-superstability of these groups in many cases. 
Non-abelian free groups were first shown to be not superstable by Gibone in his doctoral 
dissertation under the supervision of Sabbagh. We did not manage to see a copy of 
the thesis (in principle kept somewhere in Paris), but we guess that 
it was a rather involved argument on trees of formulas. 
Then Poizat gave in \cite{MR704088}
an argument based on genericity, and later an argument by linearity in dimension $2$ has been given 
in \cite{MR2206852}. One advantage of the argument by genericity in  \cite{MR704088} is that it 
works for all free products of two non-trivial groups (with the single exception of the dihedral case). 
Our arguments here for generic types are adaptations of those in \cite{MR704088} 
concerning free products. 

In a second type of results, we prove that most free products of groups without amalgamation 
are connected, i.e., 
with no proper definable subgroups of finite index. Again this will be a mere adaptation 
of an argument contained in \cite{MR704088}, concerning the free group on countably 
many generators. Together with a deep elementary 
equivalence result obtained in \cite{Sela10} about free products, we will get the following. 

\begin{theorem*}[\ref{TheoFreeProdDefCon}]
Let $G$ and $H$ be two nontrivial groups. Then a group elementarily equivalent to 
$G\ast H$ is connected except when $G\simeq H\simeq {\Z / 2\Z}$.
\end{theorem*}

We thank Bruno Poizat for his patient explainations of the argument for non-superstability 
in \cite[\S7]{MR704088}, and Chlo\'e Perin for her vigilance on a widely erroneous previous version 
of this note. We also thank Abderezak Ould Houcine and the anonymous referee 
for helpful comments. 

\section{Generic types in generic cases}

As in stable group theory, we say that a subset $X$ of a group $G$ is {\em left-generic} 
in $G$ when finitely many left-translates of $X$ cover $G$; {\em right-generic} 
sets are defined likewise. Notice that $G$-invariant sets $X$ have no political 
opinion: $gX=Xg$ for any $g$, so that left, right, or even bipartisan genericity, are the 
same for $X$. 

The argument for the non-superstability of free products of groups in 
\cite[\S7]{MR704088} roughly splits into two parts, first a transfer from 
algebraic properties to generic properties, and then a contradiction to superstability with 
the generic properties. The following essentially extracts the driving mechanism for the 
second part. 

\begin{Proposition}\label{TheoPoizatPart2}
Let $G$ be a group, $E$ a subset of $G$, and $f$ a parameter-free definable map from 
$G$ to $G$ such that $E\cup f(G)$ is not left-generic. 
\begin{itemize}
\item[(1)] 
If $G$ is stable and $g$ is generic over $G$, 
then $g$ satisfies the formula $(\forall x ~ g\neq f(x))$. 
\item[(2)] 
Suppose there exists a finite $s$ such that 
$|f^{-1}(f(g))|\leq s$ for every $g$ in $G\setminus{(E\cup f(G))}$. Then the 
complement of the definable subset 
$Z=\{g\in G~|~{(\forall x ~ g\neq f(x))}~ \wedge~ |f^{-1}(f(g))|\leq s\}$
is not left-generic. 
In particular, if $G$ is stable and $g$ is generic over $G$, then $g$ satisfies the formula 
defining $Z$. 
\item[(3)] 
In case $(2)$, $G$ cannot be superstable. 
\end{itemize}
\end{Proposition}
\proof
(1). The complement of $\{g\in G~|~\forall x ~ g\neq f(x)\}$ is $f(G)$, and therefore 
is not left-generic. Hence the claim follows in the stable case (see \cite[\S2]{MR704088} for generic types and 
definable sets in the stable case). 

(2). 
The complement of $Z$ is in $E\cup f(G)$, and therefore is not left-generic. 
Hence the claim in the stable case follows. 

(3). Suppose $G$ stable and consider $g$ generic over $G$. By (1), $f(g)$ is not generic 
over $G$. By (2), $|f^{-1}(f(g))|\leq s$, and in particular $g$ is algebraic over 
$f(g)$. Now $G$ cannot be superstable by the weak regularity of the generic types of a 
superstable group (see \cite[p.346]{MR704088}). 
\qed

\medskip
We say that a free product of groups $G\ast_A H$ 
is {\em non-trivial} when the amalgamated subgroup $A$ is a proper subgroup of both 
factors $G$ and $H$. Following the terminology in use in Bass-Serre theory, 
elements conjugated to an element of a factor are called 
{\em elliptic}, and {\em hyperbolic} otherwise. We say that a non-trivial free product 
$G\ast_A H$ is of {\em (f,s)-type} when it satisfies the assumptions in Proposition \ref{TheoPoizatPart2}(2), 
with $E$ the set of elliptic elements. We say that 
it is of {\em generic type} when it is of $(f,s)$-type for some parameter-free definable function $f$ and 
some finite $s$. 
With this terminology Proposition \ref{TheoPoizatPart2} rephrases as follows. 

\begin{Corollary}\label{CorGENGEN}
A nontrivial free product $P=G\ast_A H$ of generic type is not superstable. If it is of 
$(f,s)$-type and stable, and $x$ is generic over $P$, then 
$f(x)$ is not generic and $|f^{-1}(f(x))|\leq s$. 
\end{Corollary}

In \cite[\S7]{MR704088}, the non-superstability of a 
non-trivial free product without amalgamation $G\ast H$, where one of the factors contains 
at least three elements, consists first in showing that it is of $(f,1)$-type, where 
$f(x)=x^2$ (so that, assuming stability, the generic is then not a square and the unique 
square root of its square). A slight modification of the argument given in 
\cite[\S7]{MR704088} implies that any such free product is indeed of $(f,1)$-type 
for any function $f$ consisting in taking $d$-th powers, with $d>1$. 
The most natural generalization in free products with amalgamation is the next lemma. 

Recall from \cite[p.187]{LyndonSchupp77} that in a free product $P=G\ast_A H$ 
with amalgamated subgroup $A$, every element can be written in {\em reduced} form, i.e., 
as an alternating product of elements successively in different factors. The {\em length} 
$|x|$ of an element $x$ of $P$  is then the common length of its reduced forms (here we adopt the 
convention that $|x|=0$ when $x\in A$ and $|x|=1$ when $x\in {(G\cup H)\setminus A}$). 
In the sequel, we say that a free product $G\ast_A H$ is {\em non-degenerate} 
when it is non-trivial and one of the factors $G$ or $H$ contains at least two 
double cosets $AxA$ of $A$. 

\begin{Lemma}\label{Lem1}
Let $P=G\ast_A H$ be a non-degenerate free product of groups and let $E$ denote 
the set of its elliptic elements. 
Fix an integer $d>1$ and, for $n\geq 0$, let $B_n=\{y\in P~:~|y|\leq n\}$. 
\begin{itemize}
 \item[(1)] For any $n\geq 0$, there exists $\alpha_n$ in $P$ such that, for every $x\in B_n$, 
 $x\alpha_n \notin {E\cup {\{y^d~|~y\in P\}}}$. 
 \item[(2)] For any $n\geq 0$, $B_n{(E\cup {\{y^d~|~y\in P\}})}$ is a proper subset of $P$. 
 In particular $E\cup {\{y^d~|~y\in P\}}$ is not generic in $P$. 
\end{itemize}
\end{Lemma}
\proof
(1). 
By symmetry, we may assume that $H$ contains two distinct double cosets of $A$. 
For the sake of completeness, let us first recall the element $\alpha_n$ given for $d=2$ when $A=1$ in 
\cite[\S7]{MR704088}. Choose 
$g\in G\setminus A$ and $h$, $h' \in H\setminus A$ with $h\neq h'$ or $h\neq {h'}^{-1}$,  
and let $\alpha_n=(gh)^{n+4}(gh'g{h'}^{-1})^{3n+3}$. 

In general, one can proceed as follows. 
Let $g\in G\setminus A$ and $h$, $h'\in H\setminus A$ such that $AhA\neq Ah'A$. 
Consider then an element $\alpha_n$ of the form $(gh)^\alpha(gh'g{h'}^{-1})^{\beta}$, 
with $\max(d,n)\ll \alpha \ll \beta$. 
For every element $x$ of length at most $n$, consider $x\alpha_n$. Operating simplifications between $x$ and 
$(gh)^\alpha$, $x\alpha_n$ can be written in reduced form as a concatenation $x'(gh)^{\alpha'}(gh'g{h'}^{-1})^{\beta}$, with 
$\max(d,|x'|) \ll \alpha' \ll \beta$ since $|x'|\leq n$ and $\alpha' \geq \alpha-n-1$. 
Recall now that a reduced form is {\em cyclically reduced} when the extremal letters are not in the 
same factor (or the length is at most one). One can conjugate the previous reduced form to 
a cyclically reduced form, by operating at most $n+1$ conjugations by 
$(gh'g{h'}^{-1})^{-1}$. For, notice that if one gets a reduced form of the 
form $g_0h(gh)(gh)\cdots (gh'g{h'}^{-1})gh'^{-1}g$, if $gg_0\in A$ then $h'^{-1}gg_0h \notin A$ since 
$AhA\neq Ah'A$. Hence $x\alpha_n$ is conjugated to a cyclically reduced concatenation $\gamma$ 
of the form $x''(gh)^{\alpha''}(gh'g{h'}^{-1})^{\beta''}$, where one can impose 
$\alpha''$ and $\beta''$ to be as big as one wants compared to $n$ and $d$. 
It suffices now to show that $\gamma$ is not in $E$ and not a $d$-th power. 

Since $\gamma$ is cyclically reduced and $|\gamma|>1$, $\gamma$ is of minimal length in 
its conjugacy class. By the Conjugacy Theorem for Free Products with Amalgamation 
\cite[p.212]{MR0207802}, every cyclically reduced conjugate of $\gamma$ is $A$-conjugated to a 
cyclic permutation of $\gamma$. Since $\gamma$ is of length $>1$, it cannot be conjugated to an element 
of length $0$ or $1$, and thus is not elliptic. Suppose now $\gamma=v^d$ for some element $v$. After conjugacy, 
we may assume $v$ cyclically reduced, so that $v^d$ appears in cyclically reduced form, exactly by 
concatenating the cyclically reduced form of $v$. By the 
Conjugacy Theorem for Free Products with Amalgamation applied with $\gamma$ and $v^d$, we get that $\gamma$ and 
$v^d$ are cyclic permutations, up to $A$-conjugacy. One 
can then identify these two reduced forms by the rules descibed in \cite[p.38]{MR566274}: corresponding 
elements should be equal modulo double cosets of $A$. 
Having chosen $\alpha$ and $\beta$ such that $\alpha''$ and $\beta''$ are sufficiently divisible by $n$ and $d$, 
we then get, if $v_H$ denotes any element of $H\setminus A$ appearing in the cyclically reduced form of $v$, that 
$AhA=Av_HA=Ah'A$ ($=A{h'}^{-1}A$ also). This is a contradiction. 

(2). In the first claim, $\alpha_n \notin x^{-1}(E\cup {\{y^d~|~y\in P\}})$. 
Since $|x|=|x^{-1}|$, this shows that 
$\alpha_n \notin B_n(E\cup {\{y^d~|~y\in P\}})$, and in particular the latter is a proper 
subset of $P$. 

As the length of finitely many elements of $P$ is uniformely bounded, it follows in particular that 
$P$ cannot be covered by finitely many left-translates of $(E\cup {\{y^d~|~y\in P\}})$. 
Notice also that $E\cup {\{y^d~|~y\in P\}}$ is $G$-invariant, so that there is no ambiguity between 
left and right genericity. 
\qed

\medskip
Recall from \cite{MR2225896} that a subset $X$ of a group $G$ is called 
{\em generous} when the $G$-invariant set $X^G$ is generic in $G$. 
Lemma \ref{Lem1} gives in particular the following, with a proof more general and direct 
than the one given in \cite[Lemma 2.12]{MR2400726}. 

\begin{Corollary}\label{Cor1}
Let $P=G\ast_A H$ be a non-degenerate free product of groups. 
Then $(G\cup H)$ is not generous in $P$. 
\end{Corollary}

It should be noted here that we realized (lately) that some of the results in the present section 
may overlap, or even be weaker than, some results in \cite{Houcine2008}. Corollary \ref{Cor1} 
is certainly one of them,  and we refer to that paper for an approach through the scope of 
the Bass-Serre theory of actions on trees. It was also already observed there that $(G\cup H)$ is generous in 
$G\ast_A H$ when $A$ has index two in both factors, so that Corollary \ref{Cor1} does not hold 
for all non-trivial free products. 
The reader can also find in \cite{Houcine2008} many examples of superstable groups $G\ast_A H$ when $A$ has index 
$2$ in both factors, as well as arguments similar to those in Lemma \ref{Lem1} in the broader context of groups 
acting on trees. 

When $P$ is stable in Lemma \ref{Lem1}, Proposition \ref{TheoPoizatPart2}(1) applied with $E=\emptyset$ and 
each function $f(x)=x^d$ ($d>1$) gives the following information about generic types. 

\begin{Corollary}\label{CorGenTypeOverFreeProd}
Let $P=G\ast_A H$ be a non-degenerate free product of groups. 
If $P$ has a stable theory and 
$x$ is generic over $P$, then $x$ is not a proper power.  
\end{Corollary}

Getting the assumptions in Proposition \ref{TheoPoizatPart2}(2) in general seems to be 
more tricky. 

\begin{Corollary}\label{CorGenTypeOverFreeProd1}
Let $P=G\ast_A H$ be a non-degenerate free product of groups of $(x \mapsto x^d,s)$-type 
for some finite $s\geq 1$ and $d>1$. 
If $P$ has a stable theory and 
$x$ is generic over $P$, then $x$ is not a proper power and the number of $d$-th roots of 
$x^d$ is at most $s$
\end{Corollary}

Before mentioning groups as in Corollary \ref{CorGenTypeOverFreeProd1}, let us first see that 
a group can be of $(f,s)$-type for certain couples $(f,s)$ but not for others. 

\begin{Example}
Let $P=(\Z \times A)\ast_A (A \times \Z)$ where $A$ is an infinite elementary abelian $2$-group, so that $P$ is also 
isomorphic to $A\times F_2$. One sees that $P$ is of $(x \mapsto x^3,1)$-type, but not of 
$(x \mapsto x^2,s)$-type for any $s$. Varying the isomorphism type of $A$ in similar free products, 
one gets all possibilities as far as $(x \mapsto x^d,s)$-types are concerned. 
\end{Example}

Examples of free products $P$ satisfying an assumption as in 
Corollary \ref{CorGenTypeOverFreeProd1} cover all cases as follows: 
\begin{itemize}
\item[(1)]
{\em $A$ is finite}. To see this, let $x$ be an hyperbolic element, which we may assume to be in cyclically reduced 
form after conjugacy and of length $>1$. Then a proper power $x^d$ of $x$ is also in cyclically reduced form, 
in particular hyperbolic, and all $d$-th roots of $x^d$ are in $C_P(x^d)$. By the structure of centralizers of non-trivial 
elements from \cite[Theorem 1(i)]{KarrassSolitar77}, the only possibility is that $C_P(x^d)$ is, up to conjugacy, 
an $HNN$-extension $\<A',t~|~{A'}^t=A'\>$ for some subgroup $A'\leq A$. In other words it is a semidirect 
product $A'\rtimes \Z$ and since $A'\leq A$ is finite, on checks easily that $C_P(x^d)$ contains a finite number of 
$d$-th roots of $x^d$. Hence $P$ is of $(x \mapsto x^d,s)$-type for any $d>1$ (and $s$ may depend on $d$, but we may 
also take a uniform $s$ when $A$ is finite). 
\item[(2)]
{\em $A$ is malnormal in $G$ (or $H$)}. We may argue as above. Now we see, using reduced forms, 
that $A'=C_A(x^d)=1$, so that $C_P(x^d)\simeq \Z$ and the single $d$-th root of $x^d$ is $x$. 
Hence we may even take $s=1$ for any $d>1$ in this case. 
\item[(3)]
In fact, the above arguments show that $G\ast_A H$ is of $(x \mapsto x^d,s)$-type ($d>1$, $s\geq 1$) in all 
cases where $|C_A(x^d)|\leq s$ when $x$ varies in the set of hyperbolic elements. 
\end{itemize}

The question of the existence of a superstable non-trivial free product 
$P=G\ast_A H$ with $A$ of index at least $3$ in one of the factors remains a very 
intriguing question. At least, it cannot be of generic type by Corollary \ref{CorGENGEN}. 
The possibilities $A$ finite or malnormal were already excluded in \cite{Houcine2008} by the more 
general study of groups acting ``acylindrically" on a simplicial tree, a generalization of malnormality specific to the 
geometric approach; we refer to that paper for more on superstable groups acting on trees. 

Another inconvenience of the algebraic approach is that we have to consider 
$HNN$-extensions separately. Let us review this second case briefly. 
Consider thus an $HNN$-extension $G^*=\<G,t~|~A^t=B \>$, and let $E$ be the set of elements in $G^*$ 
conjugated to an element of $G$. We say that 
$G^*$ is of {\em (f,s)-type} if there is a parameter-free definable function $f$ on $G^*$ and a finite 
$s\geq 1$ such that $|f^{-1}(f(g))|\leq s$ for every $g$ in $G^*\setminus{(E\cup f(G^*))}$. We say that 
$G^*$ is of {\em generic type} if it is of $(f,s)$-type for some parameter-free definable function $f$ and 
some finite $s$. The analog of Corollary \ref{CorGENGEN} is then the following. 

\begin{Corollary}\label{CorGENGEN-HNN}
An $HNN$-extension $G^*=\<G,t~|~A^t=B \>$ of generic type is not superstable. If it is of 
$(f,s)$-type and stable, and $x$ is generic over $G^*$, then 
$f(x)$ is not generic and $|f^{-1}(f(x))|\leq s$. 
\end{Corollary}

As for free products, any element in an $HNN$-extension $G^*=\<G,t~|~A^t=B \>$ can be written in reduced form, 
i.e., with no subword induced by the conjugacy relation $A^t=B$. 
The {\em length} $|x|$ of an element $x$ is defined as the common number of occurences 
of $t^{\pm 1}$ in its reduced forms \cite[Chap IV 2.1-2.5]{LyndonSchupp77}. 
In the sequel, we say that  the $HNN$-extension $G^*=\<G,t~|~A^t=B \>$ is {\em non-ascending} 
when $A\cup B$ is a proper subset of $G$; notice that this is equivalent to $A<G$ and $B<G$ (so that this 
terminology coincides with the one in use in \cite{MR2486799}). 
The following is an analog of Lemma \ref{Lem1}. 
 
\begin{Lemma}\label{Lem1-HNN}
Let $G^*=\<G,t~|~A^t=B \>$ be a non-ascending $HNN$-extension, 
and let $E$ be the union of $G^*$-conjugates of $G$. 
Fix an integer $d>1$ and, for $n\geq 0$, let $B_n=\{y\in G^*~:~|y|\leq n\}$. 
\begin{itemize}
\item[(1)] 
For any $n\geq 0$, there exists $\alpha_n$ in $G^*$ such that, for every $x\in B_n$, 
 $x\alpha_n \notin {E\cup {\{y^d~|~y\in G^*\}}}$. 
 \item[(2)] 
For any $n\geq 0$, $B_n{(E\cup {\{y^d~|~y\in G^*\}})}$ is a proper subset of $G^*$. 
 In particular $E\cup {\{y^d~|~y\in G^*\}}$ is not generic in $G^*$. 
\end{itemize}
\end{Lemma}
\proof 
(1). 
Let $g$ be an element in $G\setminus (A\cup B)$. 
One can argue as in the proof of Lemma \ref{Lem1}, by considering here an element 
$\alpha_n$ of the form $(gt)^{\alpha}(gtgt^{-1})^{\beta}$ with $\beta \gg \alpha \gg \max(d,n)$. 
One can check similarly, using the alternance of $t$ and $t^{-1}$ on the right of $\alpha_n$, 
that $x\alpha_n$ is conjugated to a 
cyclically reduced element $\gamma$ of the form $x''(gt)^{\alpha''}(gtgt^{-1})^{\beta''}$, with $\beta'' \gg \alpha'' \gg \max(d,n)$. 
Here cyclically reduced forms are those reduced forms which remain reduced up to cyclic permutation. 
There is a similar Conjugacy Theorem for $HNN$-Extensions \cite[Chap IV 2.5]{LyndonSchupp77}. 
It yields similarly that the element $\gamma$, of minimal length in its conjugacy class, cannot be conjugated to an 
element of length one, and hence is not elliptic. To see that $\gamma$ cannot be conjugated to the $d$-th power of a 
cyclically reduced element $v$, one can also argue similarly. The Conjugacy Theorem yields that a cyclic 
permutation of $v^d$ is $(A\cup B)$-conjugated to $\gamma$; in this situation 
the exponents $\pm 1$ of the corresponding occurences of $t$ must be the same \cite[p.39]{MR566274}, 
and thus here we may just use the alternance of $t$ and $t^{-1}$ on the right of $\gamma$ to get a contradiction. 

(2). Follows as item (2) in Lemma \ref{Lem1}. 
\qed

\medskip
With Lemma \ref{Lem1-HNN} we get analogs of Corollaries \ref{Cor1}-\ref{CorGenTypeOverFreeProd1} 
in the same way. 

\begin{Corollary}
Let $G^*=\<G,t~|~A^t=B \>$ be a non-ascending $HNN$-extension. 
Then $G$ is not generous in $G^*$. 
\end{Corollary}

\begin{Corollary}\label{CorGenTypeOverFreeProd-2}
Let $G^*=\<G,t~|~A^t=B \>$ be a non-ascending $HNN$-extension. 
If $G^*$ has a stable theory and 
$x$ is generic over $G^*$, then $x$ is not a proper power.  
\end{Corollary}

\begin{Corollary}\label{TheoSuperstable2}
Let $G^*=\<G,t~|~A^t=B \>$ be a non-ascending $HNN$-extension of $(x \mapsto x^d,s)$-type 
for some finite $s\geq 1$ and $d>1$. 
If $G^*$ has a stable theory and 
$x$ is generic over $G^*$, then $x$ is not a proper power and the number of $d$-th roots of 
$x^d$ is at most $s$

\end{Corollary}

Examples of $HNN$-extensions satisfying the second assumption in Corollary \ref{TheoSuperstable2} 
include any $HNN$-extension $\<G,t~|~A^t=B \>$ such that $C_G(t)$ is 
finite, or such that $A^g\cap B=1$ for any $g\in G$. 
Indeed, by the structure of centralizers in the case of 
$HNN$-extensions \cite[Theorem 1(ii)]{KarrassSolitar77}, we may argue exactly as in the case of free products with amalgamation. 

Since any expansion of a non-superstable theory is also non-superstable, the non-superstability results 
in the present section are not sentive to extra structure expanding the group language. 
Our main result in the next section will work only for pure groups a priori. 

\section{Connectedness}

In this section we prove the connectedness of most free products of groups without amalgamation. 
We first elaborate a bit on an argument which stems from \cite[Lemme 6]{MR704088} (see also \cite{0799.03040} for related arguments), and which deals with automorphisms groups rather 
than with definable subsets. 

Consider  a free product without amalgamation of the form $G=A\ast F_\omega$, 
where $A$ is a group and $F_\omega=\<e_i~|~i<\omega \>$ is the free group on 
countably many generators $e_i$. 
Call a subset of $G$ {\em almost invariant} if it is setwise 
stabilized by any automorphism of $G$ fixing pointwise $A$ together with a given finite set 
of elements of $G$. In any decomposition $A\ast \<e'_i~|~i<\omega \>$ of $G$, we have that 
every finite subset of $G$ is contained in some 
subgroup $A\ast \<e'_i~|~i<n \>$ for some finite $n$. Since the subgroup of automorphisms 
of $G$ fixing $A\ast \<e'_i~|~i<n \>$ pointwise acts transitively on $\{e'_j~|~j>n\}$, 
we see that an almost 
invariant subset of $G$ contains only a finite number of the $e'_i$, or all of them except a finite number. We also see that this latter property is independent of the choice of a sequence 
$(e'_i)_{i<\omega}$ such that that $G=A\ast \<e'_i~|~i<\omega \>$. 
We can then call an almost invariant subset {\em small} in the first case, 
and {\em big} in the second case. 

\begin{Lemma}
Let $X$ be an almost invariant subset of $G$. Then $G\setminus X$ is almost invariant and the 
following holds. 
\begin{itemize}
\item[(1)]
$X$ is big if and only if $G\setminus X$ is small. 
\item[(2)]
For any element $g$ in $G$, $gX$ and $Xg$ are almost invariant. They are big whenever 
$X$ is, and small whenever $X$ is. 
\end{itemize}
\end{Lemma}
\proof
The first claim and item (1) are clear. 

For (2), we easily see that $gX$ and $Xg$ are almost invariant: add $g$ to the given finite 
set of elements such that $X$ is setwise stabilized by the group of automorphisms pointwise fixing 
$A$ plus that given set. 
Now express $g$ with the first $n$ generators of a free generating sequence 
$\{e'_i~|~i<\omega \}$ such that $G=A\ast \<e'_i~|~i<\omega \>$. 
We see that the sequence $e'_1$, ..., $e'_n$, $ge'_{n+1}$, ..., $ge'_{m+1}$, ... 
is also a free generating sequence. It follows that $gX$ is big whenever $X$ is. Similarly, 
$Xg$ is big whenever $X$ is. Now, translating by $g^{-1}$, we see that $gX$ and $Xg$ are small 
whenever $X$ is. 
\qed

\medskip
It should be clear to our reader that the intersection of two big almost invariant sets is big. 
Consequently, the union of two small almost invariant sets is small. 

\begin{Proposition}\label{Poizat}
Let $A$ be a group, $F_\omega=\<e_i~|~i<\omega \>$ the free group on countably many 
generators $e_i$, and $G=A\ast F_\omega$ the free product without amalgamation of $A$ and 
$F_\omega$. 
\begin{itemize}
\item[(1)]
Any generic (left or right) almost invariant set is big. 
\item[(2)]
$G$ has no proper almost invariant subgroup of finite index. 
\item[(3)]
$G$ is connected, i.e., with no proper definable subgroup of finite index. 
\end{itemize}
\end{Proposition}
\proof 
(1). Suppose $X$ is generic (left or right), almost invariant, and small. 
In the case of left-genericity, $G=g_1X\cup \cdots \cup g_sX$ for finitely many 
elements $g_i$, and we 
get that $G$ is small, a contradiction. In the right case we argue similarly. 

(2). Since cosets of a proper subgroup of finite index are generic (left and right) and form a 
proper partition of $G$, the claim follows from (1). 

(3). Any definable subset $X$ is almost invariant, since it is setwise fixed by all automorphisms 
fixing pointwise $A$ together with a given set of parameters in a given first-order definition of $X$. 
Hence a coset of a proper definable subgroup of finite index cannot be definable by (2). 
\qed

\medskip
With one of the results of \cite{Sela10} we get the following. 

\begin{Theorem}\label{TheoFreeProdDefCon}
Let $G$ and $H$ be two nontrivial groups. Then a group elementarily equivalent to 
$G\ast H$ is connected except when $G\simeq H\simeq {\Z / 2\Z}$.
\end{Theorem}

\proof
Since elementary equivalence preserves connectedness, we may consider directly $G\ast H$. 

If $G$ and $H$ are cyclic of order 2, then $G\ast H$ is dihedral and in particular not connected. 
If $G\ast H$ is not dihedral, then it is elementarily equivalent to $G\ast H \ast F_\omega$ 
by \cite[Theorem 7.2]{Sela10}, which is connected by Proposition \ref{Poizat}(3). 
Since elementary equivalence preserves connectedness, it follows that 
$G\ast H$ is connected. 
\qed

\medskip
We note that in the proof of Theorem \ref{TheoFreeProdDefCon} we only used the elementary equivalence 
${G\ast H}\equiv {G\ast H \ast F_\omega}$ when $G\ast H$ is not of dihedral type. It is actually expected that 
a reworking of \cite{Sela10} ``over parameters" would imply the elementary embedding ${G\ast H}\preceq {G\ast H \ast F_1}$ in this case 
(and thus elementary embeddings ${G\ast H}\preceq {G\ast H\ast F_{\kappa}} \preceq {G\ast H\ast F_{\kappa'}}$ for all 
cardinals $\kappa \leq \kappa'$). In \cite[Proposition 8.8]{MR2794551} an elementary embedding of this type is 
used in a proof of the connectedness of non-cyclic torsion-free hyperbolic groups; one could 
also argue without such an elementary embedding as is done here. 

\medskip
If $G$ is stable in Proposition \ref{Poizat}, then we see as in 
\cite[Corollary 2.7]{MR2400726} that the sequence $(e_i)_{i<\omega}$ is a Morley sequence 
of the unique generic type $p_0$ of $G$ over $\emptyset$. 
In particular, any primitive element of $F_\omega$ is a realization of the generic type of $G$ 
over $\emptyset$.  

It follows that if $G\ast H$ is a non-trivial free product not of dihedral type and stable in 
Theorem \ref{TheoFreeProdDefCon}, primitive elements of 
$F_\omega$ realize the generic type in the elementary equivalent group $G\ast H \ast F_\omega$.  
The full characterization of the set of realizations of the generic type, as in 
\cite{MR2529900} in the free group case, seems to depend on the nature of the factors $G$ and $H$; 
it is even unclear whether the generic 
type is realized in the ``standard" model $G\ast H$. Most probably, 
one could prove that the generic type is not isolated, as in \cite{MR2791345} in the free group case. 

In connection with stability also, we wish to reproduce the following comment of the 
referee around Proposition \ref{Poizat}.  
\begin{quotation}
{\it 
In an elementary equivalent group, what we can say is that the intersection of two generic definable sets is never empty; this does not mean that if the union of two definable sets is generic, then one of them should be so, as in the stable case. This last condition is equivalent to the definability of the average type of the generating sequence: the finite side of a uniform family of definable sets should be bounded (the question makes sense also for uniform families of almost invariant sets). It is very possible that Sela's argument shows the existence of such a bound: the product may be unstable for trivial and non generic reasons, but at least its generic type would behave in a stable way.
}
\end{quotation}

Theorem \ref{TheoFreeProdDefCon} proves that all non-trivial free products of groups are connected, 
with the single exception of the dihedral case. We believe that such groups are actually definably 
simple, which would follow from the following more general conjecture. 

\begin{Conjecture}\label{ConjFinale}
 Let $G\ast H$ be a free product of two groups. Then any definable subgroup of $G\ast H$ is of one of the following type: 
\begin{itemize}
 \item the full group, 
 \item conjugated to a subgroup of one of the factors $G$ or $H$, 
 \item cyclic infinite and of hyperbolic type, or 
 \item dihedral (only in case where one of the factors contains an element of order $2$). 
\end{itemize}
\end{Conjecture}

Most probably, one way to prove Conjecture \ref{ConjFinale} could be obtained by a direct generalization 
to free products of groups of the Bestvina-Feighn notion of a {\em negligeable set} in free groups, 
and by using the quantifier elimination for definable subsets of $G\ast H$ from 
\cite{Sela10}. We refer to \cite{KarlMyasnikov} for a proof in the free group case. 

\bibliographystyle{alpha}
\bibliography{biblio}

\def\cprime{$'$}
\begin{thebibliography}{MKS66}

\bibitem[dC09]{MR2486799}
Yves de~Cornulier.
\newblock Infinite conjugacy classes in groups acting on trees.
\newblock {\em Groups Geom. Dyn.}, 3(2):267--277, 2009.

\bibitem[Dye80]{MR566274}
Joan Dyer.
\newblock Separating conjugates in amalgamated free products and {HNN}
  extensions.
\newblock {\em J. Austral. Math. Soc. Ser. A}, 29(1):35--51, 1980.

\bibitem[Jal06]{MR2225896}
Eric Jaligot.
\newblock Generix never gives up.
\newblock {\em J. Symbolic Logic}, 71(2):599--610, 2006.

\bibitem[JS10]{MR2776984}
Eric Jaligot and Zlil Sela.
\newblock Makanin-{R}azborov diagrams over free products.
\newblock {\em Illinois J. Math.}, 54(1):19--68, 2010.

\bibitem[KM11]{KarlMyasnikov}
Olga Kharlampovich and Alexey Myasnikov.
\newblock Definable subsets in a hyperbolic group.
\newblock Preprint, arXiv:1111.0577v3, 2011.

\bibitem[KS77]{KarrassSolitar77}
Abraham Karrass and Donald Solitar.
\newblock Subgroups with centre in {HNN} groups.
\newblock {\em J. Austral. Math. Soc. Ser. A}, 24(3):350--361, 1977.

\bibitem[LS77]{LyndonSchupp77}
Roger Lyndon and Paul Schupp.
\newblock {\em Combinatorial group theory}.
\newblock Springer-Verlag, Berlin, 1977.
\newblock Ergebnisse der Mathematik und ihrer Grenzgebiete, Band 89.

\bibitem[MKS66]{MR0207802}
Wilhelm Magnus, Abraham Karrass, and Donald Solitar.
\newblock {\em Combinatorial group theory: {P}resentations of groups in terms
  of generators and relations}.
\newblock Interscience Publishers [John Wiley \& Sons, Inc.], New
  York-London-Sydney, 1966.

\bibitem[MP06]{MR2206852}
Yerulan Mustafin and Bruno Poizat.
\newblock Sous-groupes superstables de {${\rm SL}_2(K)$} et de {${\rm
  PSL}_2(K)$}.
\newblock {\em J. Algebra}, 297(1):155--167, 2006.

\bibitem[OH08]{Houcine2008}
Abderezak Ould~Houcine.
\newblock Superstable groups acting on trees.
\newblock Submitted, arXiv:0809.3441v1, 2008.

\bibitem[OH11]{MR2794551}
Abderezak Ould~Houcine.
\newblock Homogeneity and prime models in torsion-free hyperbolic groups.
\newblock {\em Confluentes Math.}, 3(1):121--155, 2011.

\bibitem[Pil08]{MR2400726}
Anand Pillay.
\newblock Forking in the free group.
\newblock {\em J. Inst. Math. Jussieu}, 7(2):375--389, 2008.

\bibitem[Pil09]{MR2529900}
Anand Pillay.
\newblock On genericity and weight in the free group.
\newblock {\em Proc. Amer. Math. Soc.}, 137(11):3911--3917, 2009.

\bibitem[Poi83]{MR704088}
Bruno Poizat.
\newblock Groupes stables, avec types g\'en\'eriques r\'eguliers.
\newblock {\em J. Symbolic Logic}, 48(2):339--355, 1983.

\bibitem[Poi93]{0799.03040}
Bruno Poizat.
\newblock {Is the free group stable? (Le groupe libre est-il stable?)}.
\newblock {Berlin: Humboldt-Universit\"at, Fachbereich Mathematik}, 1993.

\bibitem[Sel10]{Sela10}
Zlil Sela.
\newblock Diophantine geometry over groups {X}: The elementary theory of free
  products of groups.
\newblock preprint: http://www.ma.huji.ac.il/$\sim$zlil/, 2010.

\bibitem[Skl11]{MR2791345}
Rizos Sklinos.
\newblock On the generic type of the free group.
\newblock {\em J. Symbolic Logic}, 76(1):227--234, 2011.

\end{thebibliography}

\end{document}